\documentclass[10pt]{amsart}
\usepackage[T1]{fontenc}
\usepackage{amsmath,amssymb,lipsum,amsthm,titlesec,fancyhdr,enumerate}
\usepackage{mathrsfs}
\usepackage{tikz}
\usepackage{pgfplots}
\pgfplotsset{compat=newest}
\usepackage[toc,page]{appendix}
\usetikzlibrary{decorations.markings}
\usetikzlibrary{chains,shapes.multipart}
\usetikzlibrary{shapes,calc, arrows,}
\usetikzlibrary{automata,positioning,fit}
\usepackage{graphicx}	
\usepackage{caption}
\usepackage{subcaption}
\usepackage{float} 
\usepackage{listings}
\usepackage{hyperref}
\usepackage{newtxtext}
\usepackage[margin=1in]{geometry}
\usepackage{newtxmath}
\title{\scshape On Ordered Fuzzy Numbers Generated by Time Series}
	\author{C.~Gast} 
    \author{ R.~Traylor$^{*}$}
    \thanks{$*$:Marquette University (corresponding author) \href{mailto:rachel.traylor@marquette.edu}
    {rachel.traylor@marquette.edu}}
	\date{\today}
 
\makeatletter
\let\runauthor\@author 
\let\runtitle\@title
\makeatother

\usepackage[backend=biber,sortcites,style=numeric,giveninits=true]{biblatex}
\DeclareNameAlias{default}{family-given}

\titleformat{\section} 
  {\centering\scshape\Large}
  {\S\thesection.}
  {6pt}
  {}

\titleformat{\subsection}[runin]
  {\scshape}
  {\thesubsection.}
  {1em}
  {\addperiod}

\newtheoremstyle{solution}
{}
{}
{\color{red}}
{}
{\bfseries}
{.}
{.5em}
{}

\theoremstyle{definition}
\newtheorem{example}{Example}[section]
\newtheorem{defin}{Definition}[section]

\theoremstyle{theorem}

\newtheorem{prop}{Proposition}[section]
\theoremstyle{solution}

\renewcommand{\kappa}{\varkappa}

\newcommand{\addperiod}[1]{#1.}

\usepackage[backend=biber,sortcites,style=numeric,giveninits=true]{biblatex}
\DeclareNameAlias{default}{family-given}

\begin{document}

\maketitle

\begin{abstract}
        This paper proposes a new trapezoidal ordered fuzzy number representation of windowed time series based on the idea of a Japanese candlestick to those previously proposed in the literature.  We define and illustrate several descriptive statistics based on the information contained in the ordered fuzzy numbers. We utilize financial trading data from three automotive companies as a case study. Further expansion can be applied to any other forms of time series data to offer further insight into many other data driven situations. 

\end{abstract}

\section{Introduction}

\subsection{Fuzzy Sets and Fuzzy Numbers}
Fuzzy numbers and fuzzy sets are used widely, particularly in inference and control, but are also used in statistical applications, such as ~\cite{MB2013}. We give a brief introduction of fuzzy sets and fuzzy numbers. 

A fuzzy set is a generalization of a classical set. In classical set theory, given a set $S$ of real numbers, we define a \textit{membership function} $\mu_{S}: \mathbb{R} \to \{0,1\}$ that denotes whether a given number $x$ is a member of $S$ or not. If $\mu_{S}(x) = 1$, $x \in S$. Otherwise, $x \notin S$. 

A fuzzy set allows for partial membership, so the membership function of a fuzzy set $A$ is $\mu_{A}: \mathbb{R} \to [0,1]$. Thus, $x \in \mathbb{R}$ can be "partly in" $A$. We call $\mu_{A}(x)$ the \textit{membership level} of $x$. 

We may also consider the \textit{level sets} or \textit{$\alpha$-cuts} of a fuzzy set $A$. For a given membership level $\alpha \in [0,1]$, the level set $A_{\alpha}$ is the interval or union of intervals such that such that all $x$ in the interval or union of intervals has membership level $\alpha$. 

 A \textit{fuzzy number} is a fuzzy subset $A$ whose membership function $\mu$ satisfies the following properties~\cite{Bede2013, KaufmanGupta1985}:
\begin{enumerate}[(i)]
\item there exists $x_{0} \in \mathbb{R}$ such that $\mu(x_{0}) = 1$
\item $\mu(tx + (1-t)y) \geq \min\{\mu(x), \mu(y)\}$ for all $t \in [0,1]$, $x,y \in \mathbb{R}$\footnote {equivalently, for levels $\alpha, \alpha'$ such that $\alpha' > \alpha$, $A_{\alpha'} \subset A_{\alpha}$}
\item $\mu$ is upper semicontinuous
\item $cl(\{x \in \mathbb{R} : \mu(x) > 0\})$ is compact, where $cl(\cdot)$ denotes the closure of a set
\end{enumerate}

What this definition requires is that we must have a \textit{core}-- at least one point with membership level 1, a \textit{support} $cl(\{x \in \mathbb{R} : \mu(x) > 0\})$ that is a closed interval, and that level sets are single intervals such that greater $\alpha-$cuts are contained in lower $\alpha-$cuts.  

Dubois and Prade~\cite{DP1978} created a more restrictive class of \textit{$L-R$ fuzzy numbers}, wherein the fuzzy number's membership function is created from two "spread" functions, $L,R:[0,1] \to [0,1]$ in the following way:

 $L,R$ are two continuous, increasing functions such that $L(0) = 0 = R(0)$, $L(1) = 1 = R(1)$. Then let $a_{0}^{-} \leq a_{1}^{-} \leq a_{1}^{+} \leq a_{0}^{+}$ all be real numbers. The fuzzy set $u : \mathbb{R} \to [0,1]$ is an \textit{L-R fuzzy number} if
\[u(x) = \begin{cases}0, & x < a_{0}^{-} \\
                  L\left(\frac{x-a_{0}^{-}}{a_{1}^{-} - a_{0}^{-}}\right), & a_{0}^{-} \leq a < a_{1}^{-} \\
                1, & a_{1}^{-} \leq x < a_{1}^{+} \\
             R\left(\frac{a_{0}^{+}-x}{a_{0}^{+} - a_{1}^{+}}\right), & a_{1}^{+} \leq x < a_{0}^{+} \\
         0, & a_{0}^{+} \leq x
\end{cases}\]

The parts of an $L-R$ fuzzy number are as follows: 

\begin{itemize}
\item the \textit{core} is  $[a_{1}^{-},a_{1}^{+}]$ 
\item the \textit{support} is  $[a_{0}^{-}, a_{0}^{+}]$  
\item  the \textit{left spread}  is $\underline{a} = a_{1}^{-} - a_{0}^{-}$, and 
\item the \textit{right spread} is $\bar{a} = a_{0}^{+}- a_{1}^{+}$ . 
\end{itemize}

For example, if $L$, $R$ are linear, then we'll notice that the membership function of an $L-R$ fuzzy number yields a trapezoid, and hence we call these \textit{trapezoidal fuzzy numbers}. Specifically, the membership function of a trapezoidal fuzzy number is given by 

$$u(x) = \begin{cases}0, & x < a_{0}^{-} \\
                  \frac{x-a_{0}^{-}}{a_{1}^{-} - a_{0}^{-}}, & a_{0}^{-} \leq a < a_{1}^{-} \\
                1, & a_{1}^{-} \leq x < a_{1}^{+} \\
             \frac{a_{0}^{+}-x}{a_{0}^{+} - a_{1}^{+}}, & a_{1}^{+} \leq x < a_{0}^{+} \\
         0, & a_{0}^{+} \leq x
\end{cases}$$

\subsection{Ordered Fuzzy Numbers}
Kosinski et al~\cite{Kosinski2003, Kosinski2003b, Kosinski2006} created \textit{ordered fuzzy numbers} (OFNs) to circumvent various issues with L-R fuzzy numbers,uncontrollable results in repeatedly applied operations, which are caused by the need of intermediate approximations~\cite{Wagenknecht2001}. The principal difference between OFNs and L-R fuzzy numbers is the notion of orientation. 

An ordered fuzzy number as defined by Kosinski et al~\cite{Kosinski2003} as an ordered pair of functions $A = (\mu_{A}^{\uparrow}, \mu_{A}^{\downarrow})$, where $\mu_{A}^{\uparrow}, \mu_{A}^{\downarrow}: [0,1] \to \mathbb{R}$ and which indicate the level sets. We generally assume for applications that $\mu_{A}^{\uparrow}, \mu_{A}^{\downarrow}$ are invertible, and thus the inverses yield the piecewise membership function of the ordered fuzzy number $A$\footnote{It should be noted that Kosinski et al do not explicitly require this.}.

The orientation is a unique feature to the ordered fuzzy number, giving what can be interpreted as a time component to the ordered fuzzy number. A particle traverses the boundary of the OFN by moving "up" $\mu_{A}^{\uparrow}$ from $\alpha = 0$ to $\alpha = 1$ and "down" $\mu_{A}^{\downarrow}$ from $\alpha = 1$ to $\alpha = 0$. 

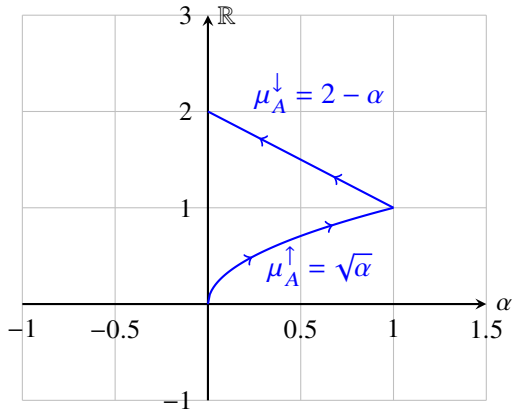
\begin{figure}[H]
\begin{tikzpicture}
\begin{axis}[scale=0.9, xmin=-1, xmax=1.5, ymin=-1, ymax=3, grid = both,
axis x line=middle,thick,
axis y line=middle,tick style={draw=none},
xlabel=$\alpha$, ylabel=$\mathbb{R}$,
every axis y label/.style={
    at={(ticklabel* cs:1.00)},
    anchor=west,
},
every axis x label/.style={
    at={(ticklabel* cs:1.00)},
    anchor=west,
}
]
\addplot[domain=0:1,samples=100,blue,
postaction={decorate, decoration={markings,
mark=at position 0.31 with {\arrow{>};},
mark=at position 0.71 with {\arrow{>};},
      }}
        ]{sqrt(x)}node[right,pos=0.9]{};
\node[below,text=blue,font=\large] at (0.6,0.7) {$\mu_{A}^{\uparrow} = \sqrt{\alpha}$};
\addplot[domain=0:1,samples=100,blue,
postaction={decorate, decoration={markings,
mark=at position 0.31 with {\arrow{<};},
mark=at position 0.71 with {\arrow{<};},
    }}
        ]{2-x}node[right,pos=0.9]{};
\node[below,text=blue,font=\large] at (0.6,2.5) {$\mu_{A}^{\downarrow} = 2-\alpha$};
\end{axis}
\end{tikzpicture}
\caption{An OFN $A = (\sqrt{\alpha}, 2-\alpha)$}
\label{fig: OFN 1}
\end{figure}
Refer to Fig~\ref{fig: OFN 1} for an example. Let $\mu_{A}^{\uparrow}(\alpha) = \sqrt{\alpha}$, so that $\mu_{A}^{\uparrow}: [0,1] \to [0,1]$, and $\mu_{A}^{\downarrow}(\alpha) = 2-\alpha$, so that $\mu_{A}^{\downarrow}:[0,1] \to [1,2]$. Then $A = (\mu_{A}^{\uparrow}, \mu_{A}^{\downarrow})$ has level sets $A_{\alpha} = [\sqrt{\alpha}, 2-\alpha]$. Then $\text{supp}(A) = [0,2]$, $A_{1/4} = [1/2, 7/4]$, and $A_{1} = \{1\}$. 

\textit{Remark: the monikers "up" and "down" do not pertain to the increasing or decreasing nature of $\mu_{A}^{\uparrow}, \mu_{A}^{\downarrow}$. They only pertain to the traversal and hence the orientation. Swapping the positions of $\mu_{A}^{\uparrow}$ and $\mu_{A}^{\downarrow}$ in the ordered pair does not change the shape, only the orientation.}

\subsection{Financial Application -- Ordered Fuzzy Candlesticks}
\label{subsec: OFC}

\begin{figure}[H]
\begin{tikzpicture}
\draw[black, thick] (-.75,13.94) -- (.75,13.94);
\draw[black, thick] (-.75,12.76) -- (.75,12.76);
\draw[black,thick] (-.75,13.94) -- (-.75,12.76); 
\draw[black,thick] (.75,13.94) -- (.75,12.76); 
\draw[black,thick] (0,14.29) -- (0,13.94);
\draw[black, thick] (0,12.76) -- (0,12.66);
\draw[black, thick] (1, 13.94) -- (1.3,13.94);
\draw[black, thick] (1,14.29) -- (1.3, 14.29);
\draw[black, thick] (1, 12.76) -- (1.3, 12.76);
\draw[black, thick] (1, 12.6) -- (1.3, 12.4);
\node at (2.5,13.94) {Close: \$1,394};
\node at (2.5, 12.76) {Open: \$1,276};
\node at (2.45, 12.3) {Low: \$1,266};
\node at (2.5, 14.29) {High: \$1,429};
\end{tikzpicture}
\caption{TSLA 7/6 - 7/10}
\label{fig: TSLA candlestick}
\end{figure}
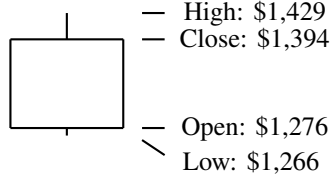

Marszalek and Burczynski~\cite{MB2013} created an OFN representation of the Japanese 
candlestick common in financial trading, as did Piasecki and Lyczkowska-Hancowiak~\cite{Piasecki2020}. 

The Japanese candlestick is a common tool used to summarized movement of time-indexed data (typically stock prices) over a given time interval (typically a single day). It resembles a boxplot, where the "box" is defined by the open and close prices of the time interval, and the "whiskers" are the maximum and minimum prices over that time interval. The box is colored red (or black/filled) if the close price is lower than the open price, and colored green (or white/unfilled) if the close price is higher than the open price. Figure~\ref{fig: TSLA candlestick} gives the Japanese candlestick for the Tesla unadjusted stock price for the week of 07/06/2020.

Given a time series $\{X_{t}\}_{t=1}^{n}$ (or a windowed subset), Marszalek and Burcynski define the \textit{ordered fuzzy candlestick} as an OFN $A = (\mu_{A}^{\uparrow}, \mu_{A}^{\downarrow})$ that is built from several parameters and a function from a chosen class. We give an overview here.

The core of the ordered fuzzy candlestick is defined by $[S_{1}, S_{2}]$, where $S_{1}$ and $S_{2}$ are chosen from the set of averages $\{\bar{X}_{SA}, \bar{X}_{LWA}, \bar{X}_{EA}\}$ such that $S_{1} < S_{2}$. $\bar{X}_{SA}$ is the simple average, $\bar{X}_{LWA}$ is a linear weighted average, and $\bar{X}_{EA}$ is an exponential average. For explicit definitions, see Appendix~\ref{sec: appendix}.  

The functions $\mu_{A}^{\uparrow}$, $\mu_{A}^{\downarrow}$ are chosen by the practitioner from a particular class to create the desired shape (trapezoidal, Gaussian, or exponential), and are a smoothed representation of the data outside $[S_{1}, S_{2}]$. For this paper, we will focus on trapezoidal OFNs, where $\mu_{A}^{\uparrow}, \mu_{A}^{\downarrow}$ are linear. 

The OFN is a long candlestick, or increasing, if $X_{1} \leq X_{n}$, and a short candlestick, or decreasing, if $X_{1} > X_{n}$. Marszalek and Burcynski determine the support as follows. 

Let $A$ be the mass above $S_{2}$ in the original time series, and $B$ the mass below $S_{1}$. Let $F^{\uparrow}$, $F^{\downarrow}$ be the areas underneath $\mu_{A}^{\uparrow}$, $\mu_{A}^{\downarrow}$, respectively. Since the functions defining the OFN are smoothed representations or estimations of the densities that enclose $A$ and $B$, then the masses $F^{\uparrow}$ and $F^{\downarrow}$ will necessarily differ from the masses they are designed to approximate. Marszalek and Burcynski in~\cite{MBorig} proposed to define the endpoints of the support of the ordered fuzzy candlestick depending on whether the OFN is long or short. 

For long ordered fuzzy candlesticks, $\mu_{A}^{\uparrow}(0) = \min\{X_{t}\} - \sigma$, where $\sigma$ is the standard deviation of the time series, and $\mu_{A}^{\downarrow}(0)$ is chosen such that 
\[\frac{F^{\uparrow}}{B} = \frac{F^{\downarrow}}{A}\]
giving the defining functions as 
\begin{equation*}
\begin{aligned}
\mu_{A}^{\uparrow}(\alpha) &= (S_{1} - [\min\{X_{t}\} - \sigma])\alpha + [\min\{X_{t}\} - \sigma] \\
\mu_{A}^{\downarrow}(\alpha) &= (S_{2} - c_{0}^{\downarrow})\alpha + c_{0}^{\downarrow}, \textit{ where } \\
c_{0}^{\downarrow} &= \frac{A}{B}(S_{1} - [\min\{X_{t}\} - \sigma]) + S_{2}
\end{aligned}
\end{equation*}

For short ordered fuzzy candlesticks, $\mu_{A}^{\uparrow}(0) = \max\{X_{t}\} + \sigma$, and $\mu_{A}^{\downarrow}(0)$ is chosen such that 
\[\frac{F^{\uparrow}}{A} = \frac{F^{\downarrow}}{B}\]
giving the defining functions\footnote{The original paper~\cite{MB2013} contains a typographical error. We have corrected it here.} as 

\begin{equation*}
\begin{aligned}
\mu_{A}^{\uparrow}(\alpha) &= (S_{2} - [\max\{X_{t}\} + \sigma])\alpha + [\max\{X_{t}\} + \sigma] \\
\mu_{A}^{\downarrow}(\alpha) &= (S_{1} - c_{0}^{\downarrow})\alpha + c_{0}^{\downarrow}, \textit{ where } \\
c_{0}^{\downarrow} &= \frac{B}{A}(S_{2} - [\max\{X_{t}\} + \sigma]) + S_{1}
\end{aligned}
\end{equation*}

In practice, we do not have a closed-form "density function" to yield $A$ and $B$. Marszalek and Burcynski~\cite{MB2013} recommend the following statistics for $A$ and $B$ (Adding a 1 to each of these sums ensures a nonzero result for $A$ and $B$.):
\[A = 1+ \sum_{i=1}^{n}X_{i}\mathbb{I}(X_{i} \geq S_{2}) \qquad B = 1 + \sum_{i=1}^{n}X_{i}\mathbb{I}(X_{i} \leq S_{1}).\]

There are several issues with this particular method of generating ordered fuzzy numbers from a time series, primarily around the arbitrariness of various choices. First, the choices of the core endpoints are completely arbitrary, with no justification or interpretation. Secondly, there is no intuitive meaning to the orientation of the OFN, since $\mu_{A}^{\uparrow}(0)$ is always fixed at an extreme datapoint $\pm$ a tolerance value, and $\mu_{A}^{\downarrow}(0)$ is designed subsequently to ensure the "error" in $F^{\downarrow}$ in estimating its portion outside the core ($A$ or $B$ depending on orientation) is equal to that of $F^{\uparrow}$. This gives no indication as to the meaning of the direction, or what the support is designed to show. Moreover, the choice to extend $\mu_{A}^{\uparrow}(0)$ by an additional $\sigma$ in the appropriate direction is not motivated or justified. In order to give some clarity and meaning to time series-generated ordered fuzzy numbers, we propose an alternative construction of ordered fuzzy numbers from time series data, and provide intuitive interpretations as well as summary statistics that can be derived from such objects. 

We also briefly mention the representation proposed by~\cite{Piasecki2020} which only allows for trapezoidal OFNs, and defines the OFN support as $[\min\{X_{1}, X_{n}\}, \max\{X_{1}, X_{n}\}]$. For financial applications, the core is defined as the interval  $[\min\{X_{o}, X_{c}\}, \max\{X_{o}, X_{c}\}]$  with endpoints of open and close prices, ($X_{o}$ and $X_{c}$, respectively). While this may be a "direct" translation of the Japanese candlestick to a trapezoidal OFN, we do not see the use of this definition for anything besides limited financial applications.

\section{A New Definition of Time-Series Generated Trapezoidal Ordered Fuzzy Numbers}
\label{sec: modification definition}

In this section, we propose an alternative method to construct ordered fuzzy numbers from time series data. This new definition better suits the original intentions of ordered fuzzy numbers, and yields an easily interpretable visual aid, as well as a mathematical and statistical object upon which we can operate arithmetically. 

Given a time series $\{X_{t}\}$, we create windows of data of size $n$, then for each window, we generate an ordered fuzzy number as follows:

The core is given by $[\min\{\bar{X}_{SA}, \bar{X}_{WA}\}, \max\{\bar{X}_{SA}, \bar{X}_{WA}\}]$, where $\bar{X}_{WA}$ is any weighted average such that the weights $\{w_{i}\}_{i=1}^{n}$ satisfy $w_{i} \leq w_{i+1}$. 

We have the following proposition:
\begin{prop}
Let $\bar{X}_{SA} = \sum_{i=1}^{n}\frac{X_{i}}{n}$, and $\bar{X}_{WA} = \sum_{i=1}^{n}\frac{w_{i}X_{i}}{\sum_{i=1}^{n}w_{i}}$ be any weighted average such that $w_{1} \leq w_{2} \leq \ldots \leq w_{n}$. If $X_{1} \leq X_{n}$, then $\bar{X}_{SA} \leq \bar{X}_{WA}$, and if $X_{1} > X_{n}$, $\bar{X}_{SA} \geq \bar{X}_{WA}$.
\label{prop: weighted average}
\end{prop}

The proof of this proposition is in Appendix~\ref{sec: appendix}.

Thus, the OFN is increasing if $X_{1} \leq X_{n}$, or, equivalently, if $\bar{X}_{SA} \leq \bar{X}_{WA}$.  The OFN is decreasing if $X_{1} > X_{n}$, or, equivalently, if $\bar{X}_{SA} > \bar{X}_{WA}$. 

Denoting $S_{1} = \min\{\bar{X}_{SA}, \bar{X}_{WA}\}$ and $S_{2} = \max\{\bar{X}_{SA}, \bar{X}_{WA}\}$, then let $I_{S_{2}}$ be the union of time intervals where $X_{t} > S_{2}$. Similarly, let $I_{S_{1}}$ be the union of intervals where $X_{t} < S_{1}$. 

We define the support of the trapezoidal OFN as $[a_{0}^{-}, a_{0}^{+}]$, where $a_{0}^{-}$ is the center of mass of $I_{S_{1}}$ and $a_{0}^{+}$ is the center of mass of $I_{S_{2}}$. For discrete time series, the support is therefore 
\[\left[\frac{\sum_{t=1}^{n}X_{t}\mathbb{I}(X_{t} < S_{1})}{\sum_{t=1}^{n}\mathbb{I}(X_{t} < S_{1})}, \frac{\sum_{t=1}^{n}X_{t}\mathbb{I}(X_{t} > S_{2})}{\sum_{t=1}^{n}\mathbb{I}(X_{t} > S_{2})}\right] \]

Thus, a long (increasing) trapezoidal OFN $A = (\mu_{A}^{\uparrow}, \mu_{A}^{\downarrow})$, where $(X_{1} \leq X_{n})$ is given by 
\[\mu_{A}^{\uparrow}(\alpha) = \alpha(S_{1} - a_{0}^{-}) + a_{0}^{-} \qquad \mu_{A}^{\downarrow}(\alpha) = a_{0}^{+} - \alpha(a_{0}^{+} - S_{2})\]
and a short (decreasing) trapezoidal OFN is given by 
\[\mu_{A}^{\uparrow}(\alpha) = a_{0}^{+} - \alpha(a_{0}^{+} - S_{2}) \qquad \mu_{A}^{\downarrow}(\alpha) = \alpha(S_{1} - a_{0}^{-}) + a_{0}^{-}\]

\begin{example}
\begin{figure}[H]
\begin{subfigure}{0.4\textwidth}
    \centering
    \includegraphics[scale=0.5]{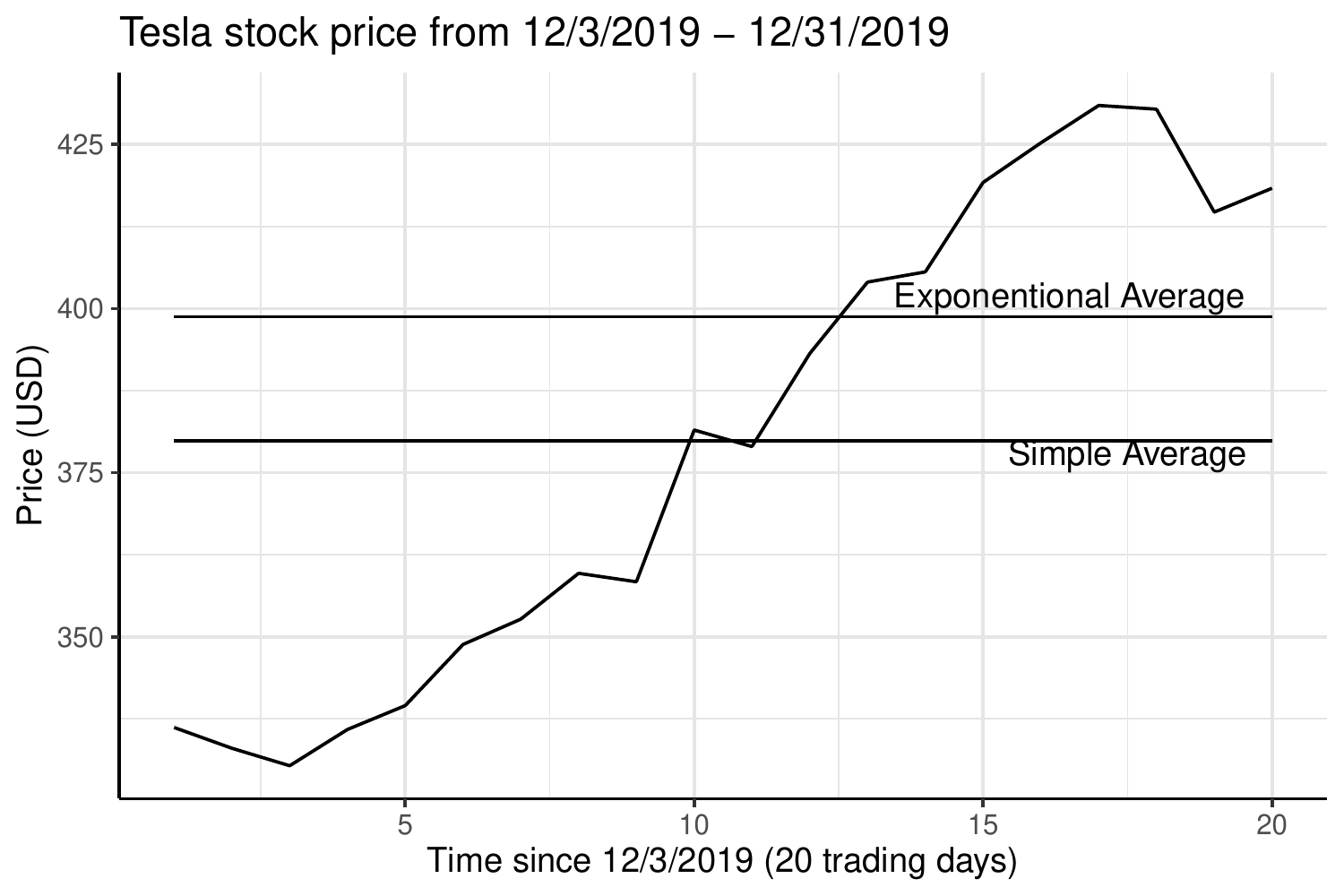}
    \caption{Tesla close prices}
    \label{fig: Tesla close}
\end{subfigure}
\hfill
\begin{subfigure}{0.4\linewidth}
\begin{tikzpicture}
\begin{axis}[scale=0.8, xmin=0, xmax=1.2, ymin=300, ymax=425, extra y ticks={325,375,425}, grid = both,
axis x line=middle,thick,
axis y line=middle,tick style={draw=none},
xlabel=$\alpha$, ylabel=$\mathbb{R}$,
every axis y label/.style={
    at={(ticklabel* cs:1.00)},
    anchor=west,
},
every axis x label/.style={
    at={(ticklabel* cs:1.00)},
    anchor=west,
}
]
\addplot[domain=0:1,samples=100,blue,
postaction={decorate, decoration={markings,
mark=at position 0.31 with {\arrow{>};},
mark=at position 0.71 with {\arrow{>};},
      }}
        ]{32.48*x + 347.36}node[right,pos=0.9]{};
\node[below,text=blue,font=\large] at (0.6,0.7) {$\mu_{A}^{\uparrow} = \sqrt{\alpha}$};
\addplot[domain=0:1,samples=100,blue,
postaction={decorate, decoration={markings,
mark=at position 0.31 with {\arrow{<};},
mark=at position 0.71 with {\arrow{<};},
    }}
        ]{418.56-19.79*x}node[right,pos=0.9]{};
\node[below,text=blue,font=\large] at (0.6,2.5) {$\mu_{A}^{\downarrow} = 2-\alpha$};
\addplot[domain=0:1,samples=100,blue,
postaction={decorate, decoration={markings,
mark=at position 0.31 with {\arrow{>};},
mark=at position 0.71 with {\arrow{>};},
    }}
        ] coordinates{(1,379.84)(1,398.77)};
\end{axis}
\end{tikzpicture}
\caption{$A_{\text{TSLA}} = (32.48\alpha + 347.36, 418.56-19.79\alpha)$}
\label{fig: tsla 2019-1}
\end{subfigure}
\caption{An example of OFN construction for the window 12/3/2019 - 12/31/2019.}
\end{figure}
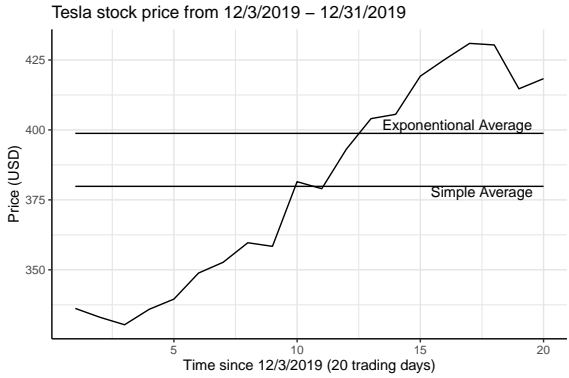
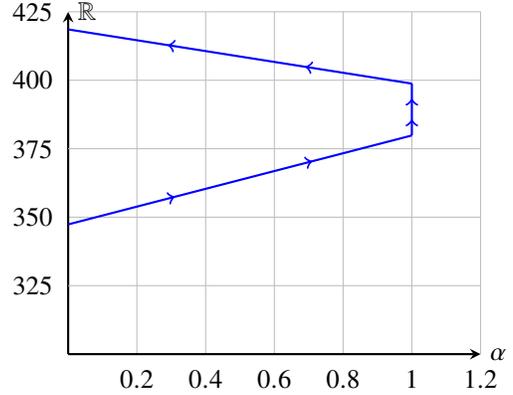

We provide an example of the construction of an OFN from time series data. Figure~\ref{fig: Tesla close} gives the daily closing prices for the trading days between 12/3/2019 and 12/31/2019. The simple average and exponential average are marked in the figure. The core is given by the exponential average and the simple average. The lower support endpoint is calculated by taking the simple average of all points below the lower endpoint of the core ($\bar{X}_{SA})$), and the upper endpoint of the support is calculated by taking the simple average above the upper endpoint of the core ($\bar{X}_{EA}$ in this case). The resulting OFN is given in Figure~\ref{fig: tsla 2019-1}. 
\end{example}

\section{Interpretation}
\label{sec: Interpretation}

The definition of the trapezoidal OFN in Section~\ref{sec: modification definition} is designed to be an estimator of typical movement over the given time interval at varying degrees of certainty $\alpha$, where $\alpha = 0$ represents the least certainty, and $\alpha = 1$ the most. In this section, we explain how to interpret the ordered fuzzy number generated from time series. 

\subsection{Core, Support, and Level Sets}
The core of a fuzzy number is classically interpreted as the interval at the highest level of membership or certainty for the concept the fuzzy number is meant to represent. For a classical fuzzy number constructed to represent the linguistic concept "a cold room", with units in degrees Fahrenheit, the core of that fuzzy number would be the interval of temperatures about which we are most certain a person would consider "cold" for a room. 

Without considering the orientation, the core of our proposed OFN is built from two averages -- a simple average and a weighted average. Both are point estimates of the center of the data in that particular time window (though both are not necessarily unbiased). In a more traditional study of time series, there would be some debate over which measure of centrality to choose, and the conclusion would be application-specific. $\bar{X}_{SA}$ assumes an equal importance or weight of all data points, regardless of time index, and $\bar{X}_{WA}$ grants greater importance or weight to the more recent points. Weighted averages towards the more recent data points reflect sudden changes in the time series more quickly than simple averages, but are also less robust and more susceptible to extreme recent values. Practitioners in finance, for example, use both types of averages, depending on the aggressiveness of their portfolios or decisions. 

Instead of giving only one point estimate for the centrality and then generating a confidence interval about that point estimate, we can give the "interval generalization" of a point estimate in the core of the generated OFN, since it is constructed from two point estimates. If desired, one can still construct confidence intervals around each endpoint of the core, giving an OFN of Type 2 found in~\cite{KaufmanGupta1985} in order to capture additional uncertainty. 

The core of this OFN also gives an indicator of movement over the period of time. By Proposition~\ref{prop: weighted average}, the orientation is determined by the ordering of $\bar{X}_{WA}$ and $\bar{X}_{SA}$. The further away $\bar{X}_{WA}$ is from $\bar{X}_{SA}$, the more movement in a particular direction. Since the core is at $\alpha = 1$, the length of the core gives an estimate of the typical movement about which we are most certain. This also helps encapsulate some of the volatility in the time series as "typical movement" or "typical volatility", so we may separate the added or additional volatility or movement. 

Generalizing the reasoning of the previous section to all level sets of an ordered fuzzy number, we may view each level set as an estimate of typical values and of movement at different levels of certainty.

Just as the core is an estimate of the typical movement about which are are most certain, the support as we have proposed in Section~\ref{sec: modification definition} should be interpreted as the widest "reasonable estimate" of typical movement at the lowest level of certainty. 

\textit{Remark: It should be noted that the support may not encompass the full range of data over the time interval. Our OFNs are designed to be a fuzzy number "point estimate" of central or typical values and movement at various levels of certainty. Thus, the support endpoints are generated by the average volatility on either side of the core endpoints in order to encapsulate volatility in a more robust manner.}

\subsection{Total Imprecision}

 Gonzales et al~\cite{Gonzalez1999} define the area below the membership function of a fuzzy number as the imprecision of the fuzzy number. This definition is sensical to use when a fuzzy number is meant to represent the uncertainty around a crisp number, but this  isn't appropriate for an ordered fuzzy number designed to estimate typical movement and the uncertainty surrounding. Given this, we propose a slight modification of the definition in~\cite{Gonzalez1999} for imprecision.

Considering the interpretations given in Section~\ref{sec: Interpretation}, we define the total imprecision of a time-series generated OFN as the total area outside the core. Formally, 

\begin{defin}[Total Imprecision]
Let $A = (\mu_{A}^{\uparrow}, \mu_{A}^{\downarrow})$ be a trapezoidal OFN constructed according to Section~\ref{sec: modification definition}. Then the \textit{total imprecision} in $A$ is given by 
\[\text{sgn}(X_{n} - X_{1})\int_{0}^{1} \mu_{A}^{\downarrow} - \mu_{A}^{\uparrow} d\alpha - [S_{2} - S_{1}]\]
 Equivalently, since the OFN is trapezoidal, the total imprecision is given by 
 \[\frac{1}{2}\left([a_{0}^{+} - a_{1}^{+}] + [a_{1}^{-} - a_{0}^{-}]\right).\]

\end{defin}

 Thus, a completely rectangular OFN would having no left or right spread, has zero imprecision.

 \textit{Remark: the notion of imprecision shouldn't be taken according to its colloquial definition. There is inherent uncertainty in even a rectangular OFN, since the core endpoints are point estimates. Here we mean imprecision as a measure of added unexplained volatility in the data.}
 
\subsection{Skew}

The next natural question is one of symmetry of the OFN. 
\begin{defin}[Skew of an OFN]
Given a trapezoidal OFN $A = (\mu_{A}^{\uparrow}, \mu_{A}^{\downarrow})$ constructed as in Section~\ref{sec: modification definition}, the \textit{skew} of $A$ is given by 
\[s = \frac{\max\{\mu_{A}^{\uparrow}(0), \mu_{A}^{\downarrow}(0)\} - S_{2}}{S_{1} - \min\{\mu_{A}^{\uparrow}(0), \mu_{A}^{\downarrow}(0)\}}\]
\end{defin}

A perfectly symmetric OFN will have $s=1$, whereas small $s$ implies greater volatility below $S_{1}$, and large $s$ implies greater volatility above $S_{2}$. There is not necessarily a time component in this metric; this is wholly derived from the shape. 

\subsection{Direction Strength}

Another useful metric from an OFN is the strength of the direction. We define the direction strength in the following way.

\begin{defin}[Direction Strength]
The \textit{direction strength} of an OFN $A = (\mu_{A}^{\uparrow}, \mu_{A}^{\downarrow})$ constructed as in Section~\ref{sec: modification definition} is given by 
\[\mathcal{D} = \frac{|\mu_{A}^{\uparrow}(1) - \mu_{A}^{\downarrow}(1)|}{\text{sgn}(X_{n} - X_{1})\int_{0}^{1}\mu_{A}^{\downarrow} - \mu_{A}^{\uparrow}d\alpha}\]
\end{defin}

It should be noted that this definition can be applied to any shape OFN, and is the percentage of total area of the OFN that belongs to the core, or highest level of certainty. 

\section{Case Study and Illustration}

To illustrate the use of OFNs in qualitatively summarizing time series data, we obtained stock price data for three companies (Tesla (TSLA), Toyota (TM), and Ford(F)) over the timeframe 10/07/2019-06/23/2020.This data includes the date, open price, close price, high and low prices, volume, and adjusted close price. All prices are measured to the cent, except for the adjusted close price, taken to 6 decimal places, and obtained from the NASDAQ website.

 Only trading dates (Monday through Friday) are included. A week here is a trading week, which is 5 days, and so forth. We illustrate the utility of data-generated OFNs using the close price over two 20 trading day periods of time (12/3/2019-13/31/2019 and 03/02/2020-03/27/2020). 

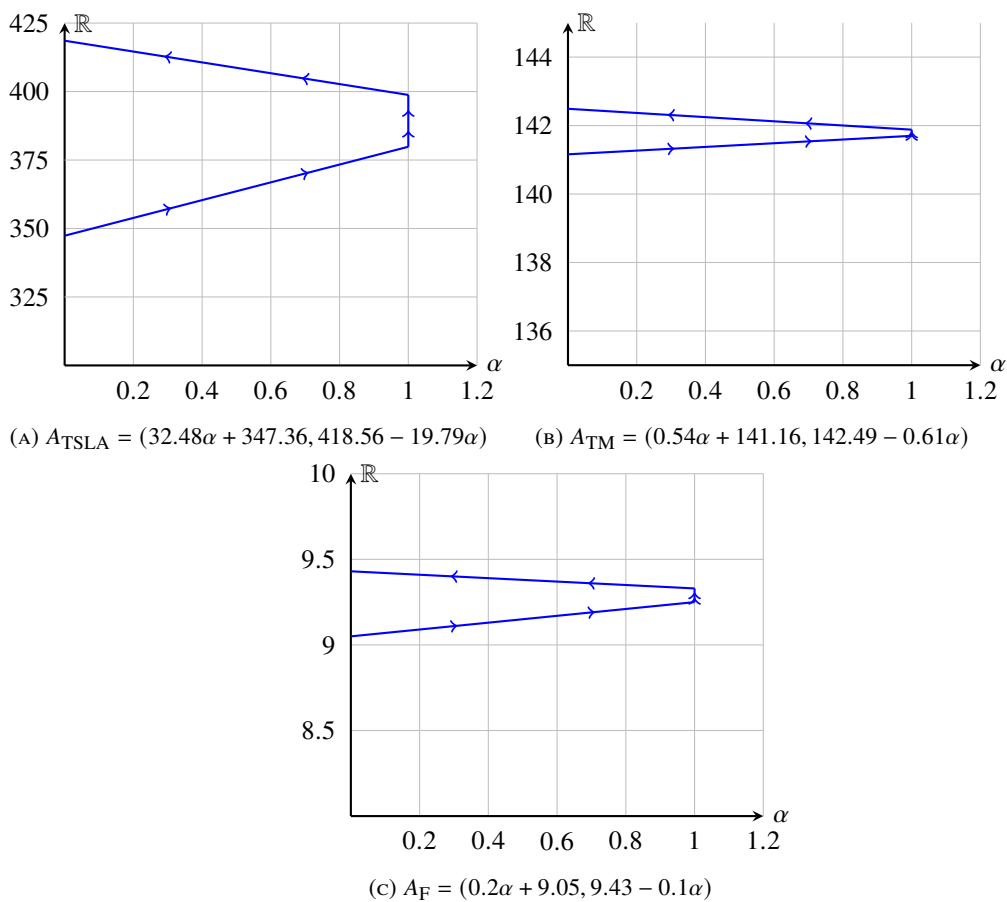
\begin{figure}[H]
\begin{subfigure}{0.4\linewidth}
\begin{tikzpicture}
\begin{axis}[scale=0.8, xmin=0, xmax=1.2, ymin=300, ymax=425, extra y ticks={325,375,425}, grid = both,
axis x line=middle,thick,
axis y line=middle,tick style={draw=none},
xlabel=$\alpha$, ylabel=$\mathbb{R}$,
every axis y label/.style={
    at={(ticklabel* cs:1.00)},
    anchor=west,
},
every axis x label/.style={
    at={(ticklabel* cs:1.00)},
    anchor=west,
}
]
\addplot[domain=0:1,samples=100,blue,
postaction={decorate, decoration={markings,
mark=at position 0.31 with {\arrow{>};},
mark=at position 0.71 with {\arrow{>};},
      }}
        ]{32.48*x + 347.36}node[right,pos=0.9]{};
\node[below,text=blue,font=\large] at (0.6,0.7) {$\mu_{A}^{\uparrow} = \sqrt{\alpha}$};
\addplot[domain=0:1,samples=100,blue,
postaction={decorate, decoration={markings,
mark=at position 0.31 with {\arrow{<};},
mark=at position 0.71 with {\arrow{<};},
    }}
        ]{418.56-19.79*x}node[right,pos=0.9]{};
\node[below,text=blue,font=\large] at (0.6,2.5) {$\mu_{A}^{\downarrow} = 2-\alpha$};
\addplot[domain=0:1,samples=100,blue,
postaction={decorate, decoration={markings,
mark=at position 0.31 with {\arrow{>};},
mark=at position 0.71 with {\arrow{>};},
    }}
        ] coordinates{(1,379.84)(1,398.77)};
\end{axis}
\end{tikzpicture}
\caption{$A_{\text{TSLA}} = (32.48\alpha + 347.36, 418.56-19.79\alpha)$}
\label{fig: tsla 2019}
\end{subfigure}
\begin{subfigure}{0.4\linewidth}
\begin{tikzpicture}
\begin{axis}[scale=0.8, xmin=0, xmax=1.2, ymin=135, ymax=145, grid = both,
axis x line=middle,thick,
axis y line=middle,tick style={draw=none},
xlabel=$\alpha$, ylabel=$\mathbb{R}$,
every axis y label/.style={
    at={(ticklabel* cs:1.00)},
    anchor=west,
},
every axis x label/.style={
    at={(ticklabel* cs:1.00)},
    anchor=west,
}
]
\addplot[domain=0:1,samples=100,blue,
postaction={decorate, decoration={markings,
mark=at position 0.31 with {\arrow{>};},
mark=at position 0.71 with {\arrow{>};},
      }}
        ]{0.54*x + 141.16}node[right,pos=0.9]{};
\node[below,text=blue,font=\large] at (0.6,0.7) {$\mu_{A}^{\uparrow} = \sqrt{\alpha}$};
\addplot[domain=0:1,samples=100,blue,
postaction={decorate, decoration={markings,
mark=at position 0.31 with {\arrow{<};},
mark=at position 0.71 with {\arrow{<};},
    }}
        ]{142.49-0.61*x}node[right,pos=0.9]{};
\node[below,text=blue,font=\large] at (0.6,2.5) {$\mu_{A}^{\downarrow} = 2-\alpha$};
\addplot[domain=0:1,samples=100,blue,
postaction={decorate, decoration={markings,
mark=at position 0.31 with {\arrow{>};},
mark=at position 0.71 with {\arrow{>};},
    }}
        ] coordinates{(1,141.70)(1,141.88)};
\end{axis}
\end{tikzpicture}
\caption{$A_{\text{TM}} = (0.54\alpha + 141.16, 142.49-0.61\alpha)$}
\label{fig: toyota 2019}
\end{subfigure}
\newline
\begin{subfigure}{0.4\linewidth}
\begin{tikzpicture}
\begin{axis}[scale=0.8, xmin=0, xmax=1.2, ymin=8, ymax=10, grid = both,
axis x line=middle,thick,
axis y line=middle,tick style={draw=none},
xlabel=$\alpha$, ylabel=$\mathbb{R}$,
every axis y label/.style={
    at={(ticklabel* cs:1.00)},
    anchor=west,
},
every axis x label/.style={
    at={(ticklabel* cs:1.00)},
    anchor=west,
}
]
\addplot[domain=0:1,samples=100,blue,
postaction={decorate, decoration={markings,
mark=at position 0.31 with {\arrow{>};},
mark=at position 0.71 with {\arrow{>};},
      }}
        ]{0.2*x + 9.05}node[right,pos=0.9]{};
\node[below,text=blue,font=\large] at (0.6,0.7) {$\mu_{A}^{\uparrow} = \sqrt{\alpha}$};
\addplot[domain=0:1,samples=100,blue,
postaction={decorate, decoration={markings,
mark=at position 0.31 with {\arrow{<};},
mark=at position 0.71 with {\arrow{<};},
    }}
        ]{9.43-0.1*x}node[right,pos=0.9]{};
\node[below,text=blue,font=\large] at (0.6,2.5) {$\mu_{A}^{\downarrow} = 2-\alpha$};
\addplot[domain=0:1,samples=100,blue,
postaction={decorate, decoration={markings,
mark=at position 0.31 with {\arrow{>};},
mark=at position 0.71 with {\arrow{>};},
    }}
        ] coordinates{(1,9.25)(1,9.33)};
\end{axis}
\end{tikzpicture}
\caption{$A_{\text{F}} = (0.2\alpha + 9.05, 9.43-0.1\alpha)$}
\label{fig: ford 2019}
\end{subfigure}
\caption{OFNs for Tesla (TSLA), Toyota (TM), and Ford (F) close prices from 12/3/2019-12/31/2019}
\label{fig: 2019 OFNs}
\end{figure}

\begin{table}[H]
\begin{tabular}{c|c|c|c|c|c|c|c|c|c|c}
Company &$a_{0}^{+}$ & $a_{0}^{-}$ & $S_{1}$ & $S_{2}$ & $\sigma$ & $X_{[1]}$ & $X_{[n]}$ & Skew & Imprecision & Direction Strength \\
\hline
Tesla & 418.57 & 347.36 & 379.84 & 398.78 & 36.68 & 336.20 & 418.33 & 0.61 & 26.13 & 0.42 \\
Toyota &142.49 & 141.16 & 141.70 & 141.89 & 0.81 & 140.73 & 140.54  &1.11 &  0.57 & 0.25 \\
Ford &  9.44 &   9.05 &   9.25 &   9.34 &  0.21 &   8.89 &  9.30 & 0.49 &   0.15 &  0.37 \\
\hline
\end{tabular}
\caption{Relevant Data from the OFNs of Tesla, Toyota, and Ford generated over 12/3/2019-12/31/2019}
\label{table: 2019 OFNs}
\end{table}

Figure~\ref{fig: 2019 OFNs} gives the data-generated OFNs for Tesla, Toyota, and Ford over the period 12/3/2019-12/31/2019, a 20 trading day period. Corresponding data is given in Table~\ref{table: 2019 OFNs}. The weighted average chosen for all three OFNs is the exponential weighted average. Since $X_{1} \leq X_{n}$ for all three stocks over this time period, all three stocks are positively oriented (increasing). Thus, $S_{1}$ is the simple average, and $S_{2}$ is the exponential average in Table~\ref{table: 2019 OFNs}. 

Tesla's stock prices are well-known for their volatility over time, and this is captured in the OFN of Figure~\ref{fig: tsla 2019}. Recall that the OFN is giving an estimate of typical movement at various levels of certainty $\alpha$. The core is given by $[379.84,398.78]$, and thus we are most certain about an approximately \$19 increase in typical stock price. The support is given by $[347.36,418.57]$. Perhaps, then, an aggressive investor would make an argument for the underlying or typical price to have increased about \$71 over the 20 day period. He would, however, have to admit that his level of certainty is at the lowest level. One can then envision the use of these OFNs for qualitative exploration or decision-making about typical price movement (with extraneous volatility removed) at varying degrees of ``aggressiveness" or "conservatism", from $\alpha = 0$ to $\alpha = 1$. 

The standard deviation of the stock price over the 20 day period is \$36.68, indicating a fairly high volatility. This is reflected in the imprecision of 26.13 given in Table~\ref{table: 2019 OFNs}, and visually we see that the OFN has a decent amount of area outside the core. The imprecision here can be interpreted as the "typical volatility" around the prices, or as the overall measure of uncertainty in the OFN itself.  

The direction strength is 0.42, indicating that the core accounts for 42\% of the total area of the OFN. We may interpret this qualitatively as a fairly strong positive orientation, despite the volatility. 

The skew is given as 0.61, which is the ratio of the area of the upper spread to the lower spread. In this case, the area below $S_{1}$ is larger than the area above $S_{2}$. One interpretation of the skew can be one of "relative volatility" in each side of the OFN. 

Similarly, we may look at the OFNs for Toyota and Ford in Figures~\ref{fig: toyota 2019} and~\ref{fig: ford 2019}. Both of these companies have much more narrow OFNs than Tesla, indicating far less volatility, as given in their respective imprecision values in Table~\ref{table: 2019 OFNs}. Both are also increasing, though not estimated to be by much. The direction strength is still somewhat comparable to the direction strength of Tesla's OFN due to the low volatility in each of Toyota's and Ford's stocks. 

A qualitative assessment of these three stocks might be that Toyota and Ford are both fairly stable stocks on a slight upward trajectory over the month of December 2019, and we are reasonably confident of the strength of that direction. Tesla has a much larger upward trajectory, but a very large imprecision to accompany it, suggesting a large amount of inherent volatility.

\begin{figure}[H]
\begin{subfigure}{0.4\linewidth}
\begin{tikzpicture}
\begin{axis}[scale=0.8, xmin=0, xmax=1.2, ymin=400, ymax=700, extra y ticks={400,450,550,650}, grid = both,
axis x line=middle,thick,
axis y line=middle,tick style={draw=none},
xlabel=$\alpha$, ylabel=$\mathbb{R}$,
every axis y label/.style={
    at={(ticklabel* cs:1.00)},
    anchor=west,
},
every axis x label/.style={
    at={(ticklabel* cs:1.00)},
    anchor=west,
}
]
\addplot[domain=0:1,samples=100,blue,
postaction={decorate, decoration={markings,
mark=at position 0.31 with {\arrow{>};},
mark=at position 0.71 with {\arrow{>};},
      }}
        ]{694.28 - 130.57*x }node[right,pos=0.9]{};
\node[below,text=blue,font=\large] at (0.6,0.7) {$\mu_{A}^{\uparrow} = \sqrt{\alpha}$};
\addplot[domain=0:1,samples=100,blue,
postaction={decorate, decoration={markings,
mark=at position 0.31 with {\arrow{<};},
mark=at position 0.71 with {\arrow{<};},
    }}
        ]{432.99 + 78.27*x}node[right,pos=0.9]{};
\node[below,text=blue,font=\large] at (0.6,2.5) {$\mu_{A}^{\downarrow} = 2-\alpha$};
\addplot[domain=0:1,samples=100,blue,
postaction={decorate, decoration={markings,
mark=at position 0.31 with {\arrow{<};},
mark=at position 0.71 with {\arrow{<};},
    }}
        ] coordinates{(1,511.26)(1,563.71)};
\end{axis}
\end{tikzpicture}
\caption{$A_{\text{TSLA}} = (694.28-130.57\alpha, 432.99+78.27\alpha)$}
\label{fig: tsla 2020}
\end{subfigure}
\begin{subfigure}{0.4\linewidth}
\begin{tikzpicture}
\begin{axis}[scale=0.8, xmin=0, xmax=1.2, ymin=100, ymax=140, grid = both,
axis x line=middle,thick,
axis y line=middle,tick style={draw=none},
xlabel=$\alpha$, ylabel=$\mathbb{R}$,
every axis y label/.style={
    at={(ticklabel* cs:1.00)},
    anchor=west,
},
every axis x label/.style={
    at={(ticklabel* cs:1.00)},
    anchor=west,
}
]
\addplot[domain=0:1,samples=100,blue,
postaction={decorate, decoration={markings,
mark=at position 0.31 with {\arrow{>};},
mark=at position 0.71 with {\arrow{>};},
      }}
        ]{127.80-6.14*x}node[right,pos=0.9]{};
\node[below,text=blue,font=\large] at (0.6,0.7) {$\mu_{A}^{\uparrow} = \sqrt{\alpha}$};
\addplot[domain=0:1,samples=100,blue,
postaction={decorate, decoration={markings,
mark=at position 0.31 with {\arrow{<};},
mark=at position 0.71 with {\arrow{<};},
    }}
        ]{114.89+4.81*x}node[right,pos=0.9]{};
\node[below,text=blue,font=\large] at (0.6,2.5) {$\mu_{A}^{\downarrow} = 2-\alpha$};
\addplot[domain=0:1,samples=100,blue,
postaction={decorate, decoration={markings,
mark=at position 0.31 with {\arrow{>};},
mark=at position 0.71 with {\arrow{>};},
    }}
        ] coordinates{(1,121.66)(1,119.70)};
\end{axis}
\end{tikzpicture}
\caption{$A_{\text{TM}} = (127.80-6.14\alpha, 114.89+4.81\alpha)$}
\label{fig: toyota 2020}
\end{subfigure}
\newline
\begin{subfigure}{0.4\linewidth}
\begin{tikzpicture}
\begin{axis}[scale=0.8, xmin=0, xmax=1.2, ymin=3, ymax=7, grid = both,
axis x line=middle,thick,
axis y line=middle,tick style={draw=none},
xlabel=$\alpha$, ylabel=$\mathbb{R}$,
every axis y label/.style={
    at={(ticklabel* cs:1.00)},
    anchor=west,
},
every axis x label/.style={
    at={(ticklabel* cs:1.00)},
    anchor=west,
}
]
\addplot[domain=0:1,samples=100,blue,
postaction={decorate, decoration={markings,
mark=at position 0.31 with {\arrow{>};},
mark=at position 0.71 with {\arrow{>};},
      }}
        ]{6.46-0.88*x}node[right,pos=0.9]{};
\node[below,text=blue,font=\large] at (0.6,0.7) {$\mu_{A}^{\uparrow} = \sqrt{\alpha}$};
\addplot[domain=0:1,samples=100,blue,
postaction={decorate, decoration={markings,
mark=at position 0.31 with {\arrow{<};},
mark=at position 0.71 with {\arrow{<};},
    }}
        ]{4.61+0.55*x}node[right,pos=0.9]{};
\node[below,text=blue,font=\large] at (0.6,2.5) {$\mu_{A}^{\downarrow} = 2-\alpha$};
\addplot[domain=0:1,samples=100,blue,
postaction={decorate, decoration={markings,
mark=at position 0.31 with {\arrow{>};},
mark=at position 0.71 with {\arrow{>};},
    }}
        ] coordinates{(1,5.58)(1,5.16)};
\end{axis}
\end{tikzpicture}
\caption{$A_{\text{F}} = (6.46-0.88\alpha, 4.61+0.55\alpha)$}
\label{fig: ford 2020}
\end{subfigure}
\caption{OFNs for Tesla (TSLA), Toyota (TM), and Ford (F) close prices from 03/02/2020-03/27/2020}
\label{fig: 2020 OFNs}
\end{figure}
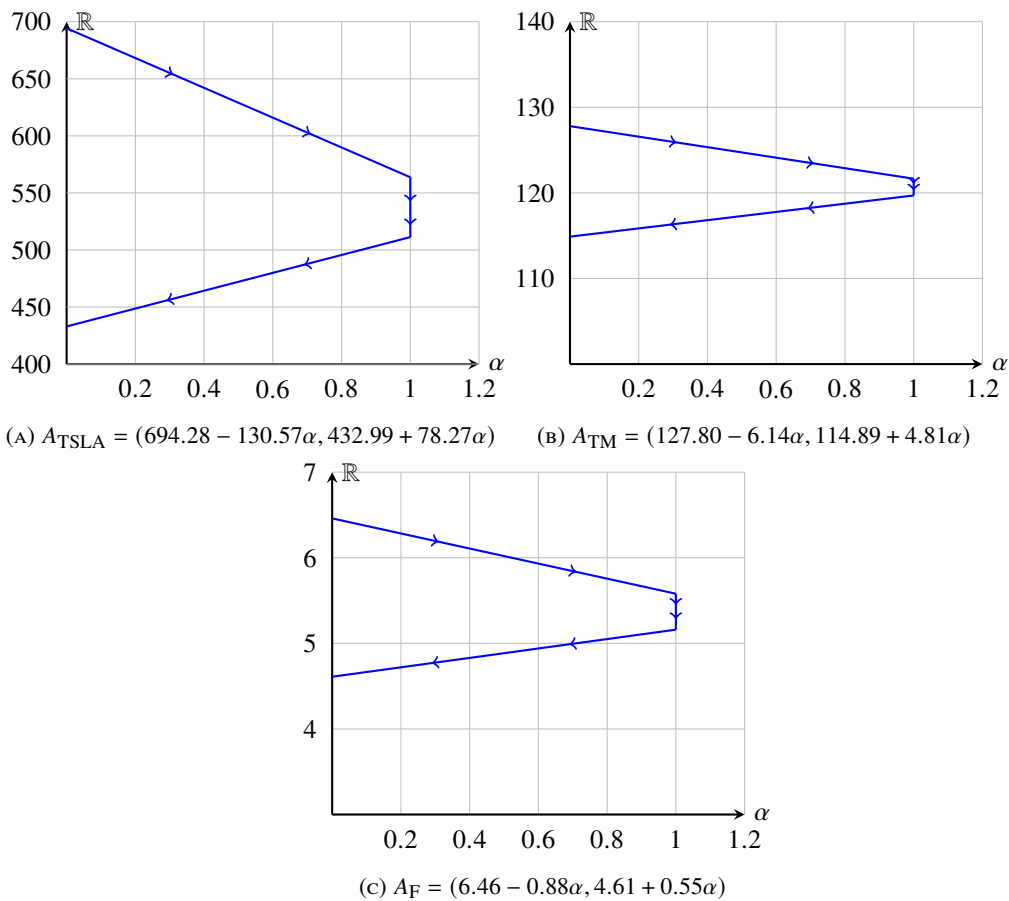

\begin{table}[H]
\begin{tabular}{c|c|c|c|c|c|c|c|c|c|c}
Company &$a_{0}^{-}$ & $a_{0}^{+}$ & $S_{1}$ & $S_{2}$ & $\sigma$ & $X_{[1]}$ & $X_{[n]}$ & Skew & Imprecision & Direction Strength \\
\hline
Tesla & 694.28 & 432.99 & 511.26 & 563.71 & 124.33 & 743.62 & 514.36 & 1.67 & 104.42 &  0.33 \\
Toyota &127.80& 114.89 & 119.70 & 121.66 &   7.23 & 132.71 & 127.24 &1.27 & 5.4721388  & 0.26 \\
Ford  & 6.46 &  4.61 &   5.16 &   5.58 &   0.96 &   7.20  & 5.19 & 1.62 &  0.7132976  & 0.37 \\
\hline
\end{tabular}
\caption{Relevant Data from the OFNs of Tesla, Toyota, and Ford generated over 03/02/2020-03/27/2020}
\label{table: 2020 OFNs}
\end{table}

Figure~\ref{fig: 2020 OFNs} and Table~\ref{table: 2020 OFNs} examine the same three stocks for a different 20 day period from 03/02/2020-03/27/2020. We can see the imprecision increase drastically for all three stock prices, which is natural given the current events of March 2020 (the beginning of the COVID-19 pandemic worldwide). All three stocks show a downward trend with reasonable direction strength. 

\section{Conclusion and Future Research}

This paper gave a new definition of time-series generated trapezoidal ordered fuzzy numbers with possible interpretations. We also defined summary metrics from the OFNs that can be used to gain qualitative information regarding the shape and movement of the time series, illustrating with a case study of three automobile manufacturer stock prices. 

Future research will include examining the statistical properties of this OFN, as well as extending this notion to other types of OFNs.

\section{Acknowledgements}

The authors have no funding sources or conflicts of interest to declare. They are exceedingly grateful to Mr. Jason Hathcock for insightful comments as well as typesetting and editing assistance. 

\begin{appendices}

\section{Appendix}
\label{sec: appendix}

\subsection{Definition of Averages}
\label{subsec: averages}
The averages used in this paper are the following.

\begin{enumerate}[(1)]
\item Simple Average:
\[\bar{X}_{SA} = \frac{\sum_{i=1}^{n}X_{i}}{n}\]
\item Linear Weighted Average (LWA):
\[\bar{X}_{LWA} = \frac{\sum_{i=1}^{n}iX_{i}}{\sum_{i=1}^{n}i} = \sum_{i=1}^{n}\frac{iX_{i}}{\left(\tfrac{n(n+1)}{2}\right)}\]
\item Exponential Average (EA): For $\gamma = \frac{2}{n+1}$,
\[\bar{X}_{EA} = \frac{\sum_{i=1}^{n}(1-\gamma)^{n-i}X_{i}}{\sum_{i=1}^{n}(1-\gamma)^{n-i}} = \sum_{i=1}^{n}\frac{\gamma(1-\gamma)^{n-i}X_{i}}{1-(1-\gamma)^{n}}\]
\end{enumerate}

\subsection{Proof of Proposition 2.1}

\textit{
 Let $\bar{X}_{SA} = \sum_{i=1}^{n}\frac{X_{i}}{n}$, and $\bar{X}_{WA} = \sum_{i=1}^{n}\frac{w_{i}X_{i}}{\sum_{i=1}^{n}w_{i}}$ be any weighted average such that $w_{1} \leq w_{2} \leq \ldots \leq w_{n}$. If $X_{1} \leq X_{n}$, then $\bar{X}_{SA} \leq \bar{X}_{WA}$, and if $X_{1} > X_{n}$, $\bar{X}_{SA} \geq \bar{X}_{WA}$.}
 
\begin{proof}

 Let the time series $\{X_{t}\}_{t=1}^{n}$, and suppose WLOG that $X_{1} \leq X_{n}$.

For $\bar{X}_{SA}$ and any weighted average, $\bar{X}_{WA}$, the expanded summations are
$$\bar{X}_{SA} = \frac{1}{n}x_1+\frac{1}{n}x_2+\frac{1}{n}x_3+\ldots + \frac{1}{n}x_n$$
$$\bar{X}_{WA} = \frac{w_1X_1}{\sum^{n}_{i=1}w_i} +\frac{w_2X_2}{\sum^{n}_{i=1}w_i}+\frac{w_3X_3}{\sum^{n}_{i=1}w_i} +\ldots+\frac{w_nX_n}{\sum^{n}_{i=1}w_i}$$

Since the sums of the weights of $\bar{X}_{SA}$ and $\bar{X}_{WA}$ are equal to each other
$$\sum_{i=1}^{n} \frac{1}{n}i = 1$$
$$\sum_{i=1}^{n}\frac{w_1}{\sum^{n}_{i=1}w_i} = 1$$

and since 
$$w_1 \leq w_2 \leq w_3 \leq \ldots \leq w_n$$

the following inequalities are true 
\begin{center}
    $\frac{w_nx_n}{\sum^{n}_{i=1}w_i} > \frac{x_n}{n}$ and $\frac{w_1x_1}{\sum^{n}_{i=1}w_i} < \frac{x_n}{n}$
\end{center}

Because of these two conditions, 
\begin{align*}
\frac{x_1}{n} &+\frac{x_2}{n}+\frac{x_3}{n}+\ldots+\frac{x_n}{n}
< \frac{w_1x_1}{\sum^{n}_{i=1}w_i}+\frac{w_2x_2}{\sum^{n}_{i=1}w_i}+\frac{w_3x_3}{\sum^{n}_{i=1}w_i}+ \ldots  +\frac{w_nx_n}{\sum^{n}_{i=1}w_i}
\end{align*}

Therefore, 
$\bar{X}_{SA} \leq \bar{X}_{WA}$. A similar argument shows the opposite inequality. 

\end{proof}
\end{appendices}
\printbibliography
\end{document}